\documentclass[12pt]{article}
\usepackage{amsmath,amsthm,amsfonts,amssymb,amscd,epsfig,euscript,version}
\makeindex

\chardef\bslash=`\\ 

\makeatletter
\def\verbatim{\interlinepenalty\@M \@verbatim
  \leftskip\@totalleftmargin\advance\leftskip2pc
  \frenchspacing\@vobeyspaces \@xverbatim}
\makeatother
\newcommand{\eul}{\EuScript}
\newcommand{\bdm}{\begin{displaymath}}
\newcommand{\edm}{\end{displaymath}}

\renewcommand{\(}{\left(}
\renewcommand{\)}{\right)}

\newcommand{\R}{{\mathbb R}}

\newcommand{\Norm}{{\eul N}}

\newcommand{\dd}[1]{\, {\mathrm d}{#1}}

\newcommand{\diag}{\operatorname{diag}}

\newcommand{\define}{\overset{\Delta}{=}}

\newcommand{\udots}{
  \mathinner {\mkern 1mu\raise 1pt \vbox {\kern 7pt \hbox {.}}\mkern 2mu
  \raise 4pt \hbox {.}\mkern 2mu\raise 7pt \hbox {.}\mkern 1mu}}

\newcommand{\bb}{{\boldsymbol b}}

\newcommand{\be}{{\boldsymbol e}}

\newcommand{\bp}{{\boldsymbol p}}

\newcommand{\br}{{\boldsymbol r}}

\newcommand{\bx}{{\boldsymbol x}}
\newcommand{\by}{{\boldsymbol y}}
\newcommand{\bz}{{\boldsymbol z}}

\newcommand{\bX}{{\boldsymbol X}}
\newcommand{\bY}{{\boldsymbol Y}}

\theoremstyle{plain}   
\newtheorem{thm}{Theorem}[section]   
\newtheorem{cor}[thm]{Corollary}     
\newtheorem{lem}[thm]{Lemma}         
\newtheorem{prop}[thm]{Proposition}  
\newtheorem{conj}[thm]{Conjecture}	

\theoremstyle{definition}
\newtheorem{rem}[thm]{Remark}        

\newcommand{\ndim}{n}

\newcommand{\nsamp}{N}

\newcommand{\Dict}{{\eul D}}

\newcommand{\train}{{\eul T}}

\newcommand{\OG}{\mathrm{O}}

\newcommand{\SL}{\mathrm{SL}}
\newcommand{\GL}{\mathrm{GL}}

\newcommand{\bone}{{\boldsymbol 1}}

\newcommand{\cond}{\,|\,}

\renewcommand{\dd}{{\mathrm d}}
\newcommand{\Cost}{{\eul C}}
\newcommand{\DS}{{\eul S}}
\newcommand{\matpart}{\, \Bigg\vert \, }
\excludeversion{hide}

\title{The Generalized Spike Process, Sparsity, and Statistical Independence}
\author{Naoki Saito\\
Department of Mathematics\\
University of California\\
Davis, CA 95616 USA\\
Email: saito@math.ucdavis.edu}
\date{}

\begin{document}
%
\maketitle
\begin{abstract}
A basis under which a given set of realizations of
a stochastic process can be represented most sparsely (the so-called best 
sparsifying basis (BSB)) and the one under which such a set becomes 
as less statistically dependent as possible
(the so-called least statistically-dependent basis (LSDB)) are 
important for data compression and have generated interests among 
computational neuroscientists as well as applied mathematicians.
Here we consider these bases for a particularly simple stochastic process
called ``generalized spike process'', which puts a single spike---whose 
amplitude is sampled from the standard normal distribution---at a random 
location in the zero vector of length $\ndim$ for each realization.

Unlike the ``simple spike process'' which we dealt with in our previous paper
and whose amplitude is constant, we need to consider the kurtosis-maximizing
basis (KMB) instead of the LSDB due to the difficulty of evaluating
differential entropy and mutual information of the generalized spike process.
By computing the marginal densities and moments, we prove that: 1)
the BSB and the KMB selects the standard basis if we restrict our basis
search within all possible orthonormal bases in $\R^\ndim$;
2) if we extend our basis search to all possible volume-preserving invertible
linear transformations, then the BSB exists and is again the standard basis
whereas the KMB does not exist.
Thus, the KMB is rather sensitive to the orthonormality of the transformations
under consideration whereas the BSB is insensitive to that.
Our results once again support the preference of the BSB over the LSDB/KMB for
data compression applications as our previous work did.
\end{abstract}
\section{Introduction}
\label{sec:intro}
This paper is a sequel to our previous paper \cite{BENICHOU-SAITO-APBOOK},
where we considered the so-called \emph{best sparsifying basis}
(BSB), and the \emph{least statistically-dependent basis} (LSDB)
for the input data which are the realizations of a very simple stochastic 
process called the ``spike process.'' 
This process, which we will refer to as the ``simple'' spike process
for convenience, puts a unit impulse (i.e., its amplitude is constant $1$)
at a random location in a zero vector of length $\ndim$.
Here, the BSB is the basis in $\R^\ndim$ that best sparsifies the given
input data, and the LSDB is the basis in $\R^\ndim$ that is the closest
to the statistically independent coordinate system (regardless of whether
such a coordinate system exists or not).  In particular, we considered
the BSB and LSDB chosen from all possible orthonormal transformations
(i.e., $\OG(\ndim)$) or all possible volume-preserving linear transformations
(i.e., $\SL^\pm(\ndim, \R)$, where any element in this set has its determinant $\pm 1$).

In this paper, we consider the BSB and LSDB for a slightly more complicated
process, the ``generalized'' spike process, and compare them with those
of the simple spike process.
The generalized spike process puts an impulse whose amplitude is sampled from 
the standard normal distribution $\Norm(0,1)$.

Our motivation to analyze the BSB and the LSDB for the generalized spike 
process stems from the work in computational neuroscience \cite{OLSHAUSEN-FIELD-NATURE}, \cite{OLSHAUSEN-FIELD-V1}, \cite{BELL-SEJNOWSKI2}, \cite{vanHateren-vanderSchaaf} as well as in computational harmonic analysis \cite{DONOHO-SCA}.
The concept of sparsity and that of statistical independence are
intrinsically different.  Sparsity emphasizes the issue of compression
directly, whereas statistical independence concerns the relationship among 
the coordinates.
Yet, for certain stochastic processes, these two are intimately related, and
often confusing.  For example, Olshausen and Field \cite{OLSHAUSEN-FIELD-NATURE}, \cite{OLSHAUSEN-FIELD-V1} emphasized the sparsity as the basis selection 
criterion, but they also assumed the statistical independence of the 
coordinates.  For a set of natural scene image patches, their algorithm
generated basis functions efficient to capture and represent
edges of various scales, orientations, and positions, which are 
similar to the receptive field profiles of the neurons in our primary visual
cortex. (Note the criticism raised by Donoho and Flesia \cite{DONOHO-FLESIA}
about the trend of referring to these functions  as ``Gabor''-like
functions; therefore, we just call them ``edge-detecting'' basis functions
in this paper.)
Bell and Sejnowski \cite{BELL-SEJNOWSKI2} used the statistical independence
criterion and obtained the basis functions similar to those of Olshausen and 
Field.
They claimed that they did not impose the sparsity explicitly and such 
sparsity \emph{emerged} by minimizing the statistical dependence among the
coordinates.
These motivated us to study these two criteria.
However, the mathematical relationship between these two criteria in the 
general case has not been understood completely.
We wish to deepen our understanding of this intricate relationship.
Therefore we chose to study such spike processes, which are much simpler
than the natural scene images viewed as a high-dimensional stochastic
process.
It is important to use simple stochastic processes first since 
we can gain insights and make precise statements in terms of theorems.
By these theorems, we now understand what are the precise conditions for
the sparsity and statistical independence criteria to select the same basis
for the spike processes, and the difference between the simple and
generalized stochastic processes.

The organization of this paper is as follows.
The next section specifies our notation and terminology.
Section~\ref{sec:spar.vs.indep} defines how to quantitatively measure 
the sparsity and statistical dependence of a stochastic process relative to 
a given basis.
Section~\ref{sec:spike} reviews the results on the simple spike process
we obtained in \cite{BENICHOU-SAITO-APBOOK}.
Our main results are presented in Section~\ref{sec:gspike} where we
deal with the generalized spike process.
We conclude with discussion in Section~\ref{sec:disc}.

\section{Notations and Terminology}
\label{sec:notation}
Let us first set our notation and the terminology.
Let $\bX \in \R^\ndim$ be a random vector with some unknown probability
density function (pdf) $f_\bX$.
Let $B \in \Dict$, where $\Dict$ is the so-called 
\index{basis dictionary} \emph{basis dictionary}.
For very high dimensional data, we often use the wavelet packets and
local Fourier bases as $\Dict$ (see \cite{SAITO-LSDB3} and references therein
for more about such basis dictionaries).
In this paper, however, we use much more larger dictionaries: 
\index{$\OG(\ndim)$} $\OG(\ndim)$
(the group of orthonormal transformations in $\R^\ndim$) 
or \index{$\SL^\pm(\ndim,\R)$} $\SL^\pm(\ndim,\R)$ (the group of invertible volume-preserving
transformations in $\R^\ndim$, i.e., their determinants are $\pm 1$).
We are interested in searching a basis under which 
the original stochastic process becomes either the sparsest or 
the least statistically dependent among the bases in $\Dict$.
\begin{hide}
Let $\Cost(B \cond \train)$ be a numerical measure of \emph{deficiency} 
or \emph{cost} of the basis $B$ given the training dataset $\train$.
For the sparsity, this cost measure the denseness of the data relative 
Under this setting,
$B_\star = \arg \min_{B \in \Dict} \Cost(B \cond \train)$ is called
the \emph{best basis} relative to the cost $\Cost$ and the training dataset
$\train$.
\end{hide}
Let $\Cost(B \cond \bX)$ be a numerical measure of \emph{deficiency} 
or \emph{cost} of the basis $B$ given the input stochastic process $\bX$.
Under this setting, \index{best basis} the \emph{best basis} for the 
stochastic process $\bX$ 
among $\Dict$ relative to the cost $\Cost$ is written as
$B_\star = \arg \min_{B \in \Dict} \Cost(B \cond \bX)$.

We also note that $\log$ in this paper implies $\log_2$, unless stated
otherwise.  The $\ndim \times \ndim$ identity matrix is denoted by
$I_\ndim$, and the $\ndim \times 1$ column vector whose entries are all
ones, i.e., $(1,1,\ldots,1)^T$, is denoted by $\bone_\ndim$.

\section{Sparsity vs. Statistical Independence}
\label{sec:spar.vs.indep}
Let us now define the measure of sparsity and that of statistical 
independence to evaluate a given basis (coordinate system).
\subsection{Sparsity}
\index{sparsity} Sparsity is a key property as a good coordinate system for compression.
The true sparsity measure for a given vector $\bx \in \R^\ndim$ is 
the so-called \index{$\ell^0$ quasi-norm} $\ell^0$ quasi-norm
which is defined as
\bdm
\| \bx \|_0 \define \#\{ i \in [1,\ndim] : x_i \neq 0 \},
\edm
i.e., the number of nonzero components in $\bx$.
This measure is, however, very unstable for even small perturbation
of the components in a vector.  Therefore, a better measure is
\index{$\ell^p$ norm} the $\ell^p$ norm:
\bdm
\label{eq:ellp-norm}
\| \bx \|_p \define \( \sum_{i=1}^\ndim |x_i|^p \)^{1/p}, \quad 0 < p \leq 1.
\edm
In fact, this is a quasi-norm for $0 < p < 1$ 
since this does not satisfy the triangle inequality, 
but only satisfies weaker conditions:
$\| \bx + \by \|_p \leq 2^{-1/p'} (\| \bx \|_p + \| \by \|_p)$ where
$p'$ is the conjugate exponent of $p$; and
$\| \bx + \by \|^p_p \leq \| \bx \|^p_p + \| \by \|^p_p$.
It is easy to show that $\lim_{p \, \downarrow \, 0} \| \bx \|_p^p = \| \bx \|_0$.
See \cite{DONOHO-SCA} for the details of the $\ell^p$ norm properties.

Thus, we can use the expected $\ell^p$ norm minimization as a criterion 
to find the best basis for a given stochastic process in terms of sparsity:
\begin{equation}
\label{eq:cp}
\Cost_{p}(B \cond \bX) = E \| B^{-1} \bX \|^p_p,
\end{equation}
We propose to use the minimization of this cost to select 
\index{best sparsifying basis (BSB)} the \emph{best sparsifying basis} (BSB):
\bdm
\label{eq:min-cp}
	B_{p}= \arg \min_{B \in \Dict} \Cost_{p}(B \cond \bX).
\edm
\begin{rem}
It should be noted that \emph{the minimization of the $\ell^p$ norm can also be
achieved for each realization}. 
Without taking the expectation in \eqref{eq:cp}, 
one can select the BSB $B_p=B_p(\bx, \Dict)$ for each 
realization $\bx$.  We can guarantee that
\bdm
\min_{B \in \Dict} \Cost_p (B \cond \bX=\bx) \leq \min_{B \in \Dict} \Cost_p (B \cond \bX) \leq \max_{B \in \Dict} \Cost_p (B \cond \bX=\bx).
\edm
For highly variable or erratic stochastic processes, however, $B_p(\bx, \Dict)$
may significantly change for each $\bx$ and we need to store more information
of this set of $\nsamp$ bases if we want to use them to compress the entire
training dataset.  Whether we should adapt a basis per realization or
on the average is still an open issue.  See \cite{SAITO-BENICHOU-LARSON-LEBORNE-LUCERO} for more details.
\end{rem}

\subsection{Statistical Independence}
\index{statistical independence}
The statistical independence of the coordinates of $\bY \in \R^\ndim$ means
$f_\bY(\by) = f_{Y_1}(y_1) f_{Y_2}(y_2) \cdots f_{Y_\ndim}(y_\ndim)$,
where $f_{Y_k}$ is a one-dimensional \index{marginal distribution}
marginal pdf of $f_\bY$.
The statistical independence
is a key property as a good coordinate system for compression and particularly
modeling because: 1) damage of one coordinate does not propagate to the others;
and 2) it allows us to model the $\ndim$-dimensional stochastic process of 
interest as a set of 1D processes.
Of course, in general, it is difficult to find a truly statistically 
independent coordinate system for a given stochastic process. 
Such a coordinate system may not even exist for a certain stochastic process.
Therefore, the next best thing we can do is to find the least-statistically
dependent coordinate system within a basis dictionary.
Naturally, then, we need to measure the ``closeness'' of a coordinate
system $Y_1,\ldots,Y_\ndim$ to the statistical independence.
This can be measured by \index{mutual information} \emph{mutual information}
or relative entropy between the true pdf $f_\bY$ and the product of its
marginal pdf's:
\begin{eqnarray*}
I(\bY)  & \define &\int f_\bY(\by) \log \frac{f_\bY(\by)}{\prod_{i=1}^\ndim f_{Y_i}(y_i)} \dd{\by} \\
        & = & - H(\bY) + \sum_{i=1}^\ndim H(Y_i),
\end{eqnarray*}
where $H(\bY)$ and $H(Y_i)$ are the \index{differential entropy} differential entropy of $\bY$ and $Y_i$
respectively:
\begin{eqnarray*}
 H(\bY) &=& - \int f_\bY(\by) \, \log f_\bY(\by) \dd{\by}\\
 H(Y_i) &=& - \int f_{Y_i}(y_i) \, \log f_{Y_i}(y_i) \dd{y_i}.
\end{eqnarray*}
We note that $I(\bY) \geq 0$, and $I(\bY) = 0$ if and only if
the components of $\bY$ are mutually independent. See \cite{COVER-THOMAS}
for more details of the mutual information.

Suppose $\bY = B^{-1} \bX$ and $B \in \GL(\ndim,\R)$ with $\det B=\pm 1$.
We denote such a set of matrices by \index{$\SL^\pm(\ndim,\R)$} $\SL^\pm(\ndim,\R)$.  
Note that the usual $\SL(\ndim,\R)$ is a subset of $\SL^\pm(\ndim,\R)$.
Then, we have
\bdm
I(\bY) = - H(\bY) + \sum_{i=1}^\ndim H(Y_i) = - H(\bX) + \sum_{i=1}^\ndim H(Y_i),
\edm
since the differential entropy is \emph{invariant}
under such an invertible volume-preserving linear transformation, i.e.,
\bdm
H(B^{-1}\bX)=H(\bX)+\log | \det B^{-1}|= H(\bX),
\edm
because $| \det B^{-1} | = 1$.
Based on this fact, we proposed the minimization of the following cost 
function as the criterion to select
the so-called \index{least statistically-dependent basis (LSDB)} \emph{least statistically-dependent basis} (LSDB) in the
basis dictionary context \cite{SAITO-LSDB3}:
\begin{equation}
\label{eq:c-lsdb}
\Cost_{H}(B \cond \bX)  = \sum_{i=1}^\ndim H\((B^{-1}\bX)_i\)
= \sum_{i=1}^\ndim H(Y_i). 
\end{equation}
Now, we can define the LSDB as
\bdm
	B_{LSDB} = \arg \min_{B \in \Dict} \Cost_{H}(B \cond \bX).
\edm
We were informed that Pham \cite{PHAM} had proposed the minimization of
the same cost \eqref{eq:c-lsdb} earlier.
We would like to point out the main difference between our 
work \cite{SAITO-LSDB3} and Pham's.  We used the basis libraries such as
wavelet packets and local Fourier bases that allow us to deal with datasets 
with large dimensions such as face images whereas Pham used more general 
dictionary $\GL(\ndim,\R)$.  In practice, however, the numerical optimization 
\eqref{eq:c-lsdb} clearly becomes more difficult in his general case
particularly if one wants to use this for high dimensional datasets.

Closely related to the LSDB is the concept of the
\index{\emph{kurtosis-maximizing basis} (KMB)}
\emph{kurtosis-maximizing basis} (KMB).
This is based on the approximation of the marginal differential entropy
\eqref{eq:c-lsdb} by higher order moments/cumulants
using the \index{Edgeworth expansion} Edgeworth expansion and was derived by
Comon \cite{COMON}:
\begin{equation}
\label{eq:kurtosis}
H(Y_i) \approx - \frac{1}{48}\kappa(Y_i) = -\frac{1}{48} ( \mu_4(Y_i) - 3 \mu^2_2(Y_i) )
\end{equation}
where $\mu_k(Y_i)$ is the $k$th central moment of $Y_i$,
and $\kappa(Y_i)$ $/$ $\mu^2_2(Y_i)$ is called the \emph{kurtosis} of $Y_i$.
See also Cardoso \cite{CARDOSO2} for a nice exposition of the various
approximations to the mutual information.
Now, the KMB is defined as follows:\footnote{Note that
there is a slight abuse of the terminology; We call the kurtosis-maximizing
basis in spite of maximizing unnormalized version (without the division by
$\mu^2_2(Y_i)$) of the kurtosis.}
\begin{equation}
\label{eq:c-kmb}
B_\kappa = \arg \min_{B \in \Dict} \Cost_\kappa(B \cond \bX)
= \arg \max_{B \in \Dict} \sum_{i=1}^\ndim \kappa(Y_i),
\end{equation}
where $\Cost_\kappa(B \cond \bX) = - \sum_{i=1}^\ndim \kappa(Y_i)$.
We note that the LSDB and the KMB are tightly related, yet can be different.
After all, \eqref{eq:kurtosis} is simply an approximation to the entropy up to
the fourth order cumulant.
We also would like to point out that Buckheit and Donoho \cite{BUCK-DONO}
independently proposed the same measure as a basis selection criterion, 
whose objective was to find a basis under which an input stochastic process
looks maximally \index{non-Gaussianity} ``non-Gaussian.''

\section{Review of Previous Results on the Simple Spike Process}
\label{sec:spike}
In this section, we briefly summarize the results of the simple spike process,
which we obtained previously. See \cite{BENICHOU-SAITO-APBOOK} for the details and proofs. 

An $\ndim$-dimensional \index{simple spike process} 
\emph{simple spike process} generates the standard
basis vectors $\{ \be_j \}_{j=1}^\ndim \subset \R^\ndim$ in a random order,
where $\be_j$ has one at the $j$th entry and all the other entries are zero.
One can view this process as a unit impulse located at a
random position between $1$ and $\ndim$.

\subsection{The Karhunen-Lo\`eve Basis}
\index{Karhunen-Lo\`eve basis} 
The Karhunen-Lo\`eve basis of this process is not unique and not useful
because of the following theorem.
\begin{prop}
\label{prop:spike-KLB}
The Karhunen-Lo\`eve basis for the simple spike process is any orthonormal
basis in $\R^\ndim$ containing the ``DC'' vector $\bone_\ndim=(1,1,\ldots,1)^T$.
\end{prop}
This theorem reminds us of \index{non-Gaussianity} non-Gaussianity of 
the simple spike process

\subsection{The Best Sparsifying Basis}
\index{best sparsifying basis (BSB)}
As for the BSB, we have the following result:
\begin{thm}
\label{thm:spike-sparse}
The BSB with any $p \in [0,1]$ for the simple spike process 
is \index{standard basis}
the standard basis
if $\Dict=\OG(\ndim)$ or $\SL^\pm(\ndim,\R)$.
\end{thm}

\subsection{Statistical Dependence and Entropy of the Simple Spike Process}
\label{sec:true-entropy}
Before considering the LSDB of this process, let us note a few specifics
about the simple spike process.
First, although the standard basis is the BSB for this process, 
it clearly does not provide the statistically independent coordinates.
The existence of a single spike at one location prohibits
spike generation at other locations. This implies that these coordinates
are highly statistically dependent.

Second, we can compute the true \index{entropy} entropy $H(\bX)$ 
for this process unlike other complicated stochastic processes.
Since the simple spike process selects one possible vector from the standard 
basis vectors of $\R^\ndim$ with uniform probability $1/\ndim$, 
the true entropy $H(\bX)$ is clearly $\log \ndim$.
This is one of the rare cases where we know the true high-dimensional
entropy of the process.

\subsection{The LSDB among $\OG(\ndim)$}
For \index{$\OG(\ndim)$} $\Dict=\OG(\ndim)$, we have the following theorem.
\begin{thm}
\label{thm:spike-OG}
The \index{least statistically-dependent basis (LSDB)} LSDB 
among $\OG(\ndim)$ is the following:
\begin{itemize}
\item for $\ndim \geq 5$, either \index{standard basis} 
the standard basis or the basis whose matrix
representation is
\begin{equation}
\label{eq:lsdb-OG5}
\frac{1}{\ndim} \left[ \begin{array}{ccccc}
 \ndim-2 &   -2   & \cdots &   -2   &  -2 \\
   -2    & \ndim-2& \ddots &        &  -2 \\
\vdots   & \ddots & \ddots & \ddots & \vdots \\
   -2    &        & \ddots & \ndim-2&  -2 \\
   -2    &   -2   & \cdots &   -2   & \ndim-2
\end{array} \right];
\end{equation}
\item for $\ndim = 4$, \index{Walsh basis} the Walsh basis, i.e.,
\bdm
\frac{1}{2}\left[ \begin{array}{cccc}
1 & 1 & 1 & 1   \\
1 & 1 & -1 & -1 \\
1 & -1 & 1 & -1 \\
1 & -1 & -1 & 1
\end{array}\right];
\edm

\item for $\ndim=3$, $\left[ \begin{array}{ccc}
\frac{1}{\sqrt 3} & \frac{1}{\sqrt 6} & \frac{1}{\sqrt 2} \\
\frac{1}{\sqrt 3} & \frac{1}{\sqrt 6} & \frac{-1}{\sqrt 2} \\
\frac{1}{\sqrt 3} & \frac{-2}{\sqrt 6} & 0
\end{array}\right]$; and

\item for $\ndim=2$, $\frac{1}{\sqrt{2}}\left[ \begin{array}{cccc}
1 & 1    \\
1 & -1
\end{array}\right]$,
and this is the only case where the true independence is
achieved.
\end{itemize}
\end{thm}

\begin{rem}
\label{rem:permute-sign}
Note that when we say the basis is a matrix as above,
we really mean that the column vectors of that matrix form the basis.
This also means that any permuted and/or sign-flipped (i.e., multiplied by
$-1$) versions of those column vectors also form the basis.  
Therefore, when we say the basis is a matrix $A$, we mean not only
$A$ but also its permuted and sign-flipped versions of $A$.
This remark also applies to all the propositions and theorems below, 
unless stated otherwise.
\end{rem}

\begin{rem}
There is an important geometric interpretation of \eqref{eq:lsdb-OG5}.
This matrix can also be written as: 
\bdm
B_{HR(\ndim)} \define I_\ndim - 2 \frac{\bone_\ndim}{\sqrt{\ndim}}\frac{\bone_\ndim^T}{\sqrt{\ndim}}.
\edm
In other words, this matrix represents the \index{Householder reflection}
\emph{Householder reflection} with
respect to the hyperplane $\{ \by \in \R^\ndim \, | \, \sum_{i=0}^\ndim y_i = 0 \}$ whose unit normal vector is $\bone_\ndim/\sqrt{\ndim}$.
\end{rem}
Below, we use the notation $B_{\OG(\ndim)}$ for the LSDB among $\OG(\ndim)$
to distinguish it from the LSDB among $\GL(\ndim,\R)$, which is denoted
by $B_{\GL(\ndim)}$.
So, for example, for $\ndim \geq 5$, 
$B_{\OG(\ndim)}=I_\ndim$ or $B_{HR(\ndim)}$.

\subsection{The LSDB among $\GL(\ndim,\R)$}
As discussed in \cite{BENICHOU-SAITO-APBOOK},
for the simple spike process, there is no important distinction
in the LSDB selection from \index{$\GL(\ndim,\R)$} $\GL(\ndim,\R)$ and 
from $\SL^\pm(\ndim,\R)$.
Therefore, we do not have to treat these two cases separately.
On the other hand, the generalized spike process in Section~\ref{sec:gspike}
requires us to treat $\SL^\pm(\ndim,\R)$ and $\GL(\ndim,\R)$ differently
due to the continuous amplitude of the generated spikes.

We now have the following curious theorem:
\begin{thm}
\label{thm:spike-GL}
\index{least statistically-dependent basis (LSDB)}
The LSDB among $\GL(\ndim,\R)$ with $\ndim > 2$ is
the following basis pair (for analysis and synthesis respectively):\\
\begin{equation}
\label{eq:spike-LSDB-SL-analysis}
B_{\GL(\ndim)}^{-1}=
\left[ \begin{array}{ccccccc}
a & a & \cdots & \cdots & \cdots & \cdots & a \\
b_2 & c_2 & b_2 & \cdots & \cdots & \cdots & b_2 \\
b_3 & b_3 & c_3 &  b_3  & \cdots  & \cdots & b_3 \\
\vdots & \vdots&  & \ddots &      &        & \vdots \\
\vdots & \vdots&  &        & \ddots &     & \vdots \\
b_{\ndim-1} & \cdots & \cdots & \cdots & b_{\ndim-1} & c_{\ndim-1} & b_{\ndim-1}\\
b_\ndim & \cdots & \cdots & \cdots & \cdots & b_\ndim & c_\ndim
\end{array} \right],
\end{equation}
\begin{equation}
\label{eq:spike-LSDB-SL-synthesis}
B_{\GL(\ndim)}=
\left[ \begin{array}{ccccc}
\(1 + \sum_{k=2}^\ndim b_kd_k \)/a & -d_2 & -d_3 & \cdots & -d_\ndim \\
-b_2 d_2/a & d_2 & 0 & \cdots & 0 \\
-b_3 d_3/a & 0 & d_3 & \ddots & \vdots \\
\vdots & \vdots & \ddots & \ddots  & 0 \\
-b_\ndim d_\ndim/a & 0 & \cdots & 0 & d_\ndim
\end{array} \right]
\end{equation}
where $a$, $b_k$, $c_k$ are arbitrary real-valued constants satisfying
$a \neq 0$, $b_k \neq c_k$, and $d_k = 1/(c_k-b_k)$, $k=2,\ldots,\ndim$.

\index{$\SL^\pm(\ndim,\R)$}
If we restrict ourselves to $\Dict=\SL^\pm(\ndim,\R)$, then
the parameter $a$ must satisfy:
\bdm
a = \pm \prod_{k=2}^\ndim (c_k - b_k)^{-1}.
\edm 
\end{thm}

\begin{rem}
\label{rem:lsdb-dense}
The LSDB such as \eqref{eq:lsdb-OG5} and the LSDB pair 
\eqref{eq:spike-LSDB-SL-analysis}, \eqref{eq:spike-LSDB-SL-synthesis}
provide us with further insight into the difference
between \index{sparsity} sparsity and \index{statistical independence}
statistical independence.
In the case of \index{Householder reflector} \eqref{eq:lsdb-OG5}, 
this is the LSDB, yet does not 
sparsify the spike process at all. In fact, these coordinates are completely
dense, i.e., $\Cost_0 = \ndim$.  We can also show that the sparsity measure
$\Cost_p$ gets worse as $\ndim \to \infty$.
More precisely, we have the following proposition.
\begin{prop}
\label{prop:spike-cp-OG}
\bdm
\lim_{\ndim \to \infty} \Cost_p\(B_{HR(\ndim)} \cond \bX \) = \left\{
\begin{array}{ll}
\infty & \quad \mbox{if $0 \leq p < 1$};\\
3 & \quad \mbox{if $p=1$}.
\end{array}
\right.
\edm
\end{prop}
It is interesting to note that this LSDB approaches to \index{standard basis}
the standard basis as $\ndim \to \infty$.  This also implies that
\bdm
\lim_{\ndim \to \infty} \Cost_p\(B_{HR(\ndim)} \cond \bX \) \neq \Cost_p\(\lim_{\ndim \to \infty} B_{HR(\ndim)} \cond \bX \).
\edm

As for the analysis LSDB \eqref{eq:spike-LSDB-SL-analysis}, 
the ability to sparsify the spike process depends on the values of $b_k$ 
and $c_k$.
Since the parameters $a$, $b_k$ and $c_k$ are arbitrary as long as 
$a \neq 0$ and $b_k \neq c_k$,
let us put $a=1$, $b_k=0$, $c_k=1$, for $k=2,\ldots,\ndim$.
Then we get the following specific LSDB pair:
\bdm
\label{eq:spike-LSDB-SL-01}
B_{\GL(\ndim)}^{-1}=
\left[ \begin{array}{cccc}
1 & 1 & \cdots & 1 \\
0 & & &  \\
\vdots & & I_{\ndim-1} & \\
0 & & & \end{array} \right],
\quad
B_{\GL(\ndim)}=
\left[ \begin{array}{cccc}
1 & -1 & \cdots & -1 \\
0 & & & \\
\vdots & & I_{\ndim-1} & \\
0 & & & \end{array} \right].
\edm
This analysis LSDB provides us with a sparse representation for the simple 
spike process (though this is clearly not better than the standard basis). 
For $\bY=B_{\GL(\ndim)}^{-1}\bX$,
\bdm
\Cost_p = E \left[ \| \bY \|_p^p \right]= \frac{1}{\ndim} \times 1 +\frac{\ndim-1}{\ndim} \times 2=2-\frac{1}{\ndim}, \quad 0 \leq p \leq 1.
\edm
Now, let us take $a=1$, $b_k=1$, $c_k=2$ for $k=2,\ldots,\ndim$ in \eqref{eq:spike-LSDB-SL-analysis} and \eqref{eq:spike-LSDB-SL-synthesis}. Then we get
\bdm
\label{eq:spike-LSDB-SL-12}
B^{-1}_{\GL(\ndim)} = \left[ \begin{array}{cccc}
1 & 1 & \cdots & 1 \\
1 & 2 & \ddots & \vdots \\
\vdots & \ddots & \ddots & 1 \\
1 & \cdots & 1 & 2 
\end{array} \right],
\quad
B_{\GL(\ndim)} = \left[ \begin{array}{cccc}
\ndim  & -1 & \cdots &  -1 \\
-1     &   &   &   \\
\vdots &   & I_{\ndim-1}  &   \\
-1     &   &   &   \\ 
\end{array} \right].
\edm
The sparsity measure of this process is:
\bdm
\Cost_p = \frac{1}{\ndim} \times \ndim +\frac{\ndim-1}{\ndim} \times \{ (\ndim-1)+2^p\} = \ndim+(2^p-1)\(1-\frac{1}{\ndim}\), \quad 0 \leq p \leq 1.
\edm
Therefore, the spike process under this analysis basis is completely dense,
i.e., $\Cost_p \geq \ndim$ for $0 \leq p \leq 1$ and the equality holds if and
only if $p=0$.
Yet this is still the LSDB.
\end{rem}

Finally, from Theorems~\ref{thm:spike-OG} and \ref{thm:spike-GL}, we can prove
the following corollary:
\begin{cor}
\label{cor:spike-GL}
There is no invertible linear transformation providing the 
statistically independent coordinates for the spike process for $\ndim > 2$.
\end{cor}

\section{The Generalized Spike Process}
\label{sec:gspike}
\begin{figure}
\centering{\epsfig{file=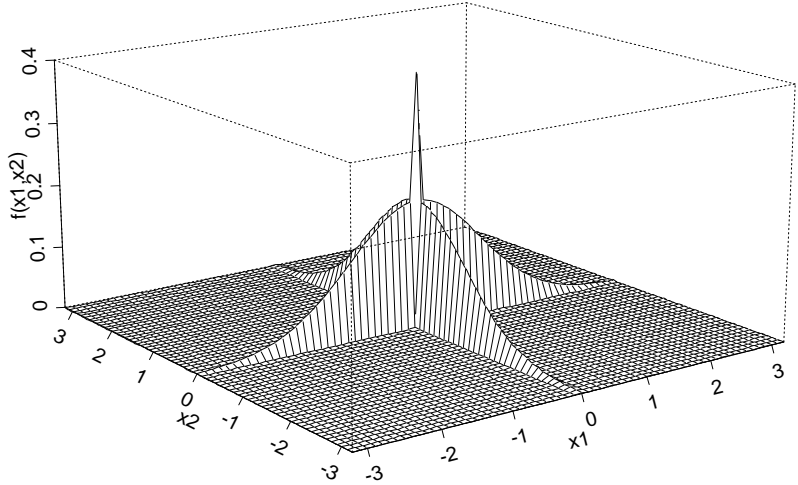,width=0.80\textwidth}}
\caption{The pdf of the generalized spike process ($\ndim=2$).}
\label{fig:genspike2}
\end{figure}

In \cite{DVDD}, Donoho et al.\ analyzed the following generalization of
the simple spike process in terms of the KLB and the rate distortion function.
This process first picks one coordinate out of $\ndim$ coordinates randomly
as before, but then the amplitude of this single spike is picked according to
the standard normal distribution $\Norm(0,1)$.
The pdf of this process can be written as follows:
\index{generalized spike process}
\begin{equation}
\label{eq:gspike}
f_\bX(\bx) = \frac{1}{\ndim} \sum_{i=1}^\ndim \( \prod_{j \neq i} \delta(x_j) \) g(x_i),
\end{equation}
where $\delta(\cdot)$ is the Dirac delta function, and $g(x)=$ $(1/\sqrt{2\pi})$
$\cdot$ $\exp(-x^2/2)$, i.e., the pdf of the standard normal distribution.
Figure~\ref{fig:genspike2} shows this pdf for $\ndim=2$.
Interestingly enough, this generalized spike process shows rather
different behavior (particularly in the statistical independence) 
from the simple spike process in Section~\ref{sec:spike}.
We also note that our proofs here are rather analytical compared to
those for the simple spike process presented in \cite{BENICHOU-SAITO-APBOOK},
which have more combinatorial flavor.

\subsection{The Karhunen-Lo\`eve Basis}
\index{Karhunen-Lo\`eve basis} 
We can easily compute the covariance matrix of this process, which is
proportional to the identity matrix.  In fact, it is just $I_\ndim/\ndim$.
Therefore, we have the following proposition, which was also stated without
proof by Donoho et al.\ \cite{DVDD}:
\begin{prop}
\label{prop:gspike-KLB}
The Karhunen-Lo\`eve basis for the generalized spike process is any 
orthonormal basis in $\R^\ndim$.
\end{prop}
\begin{proof}
Let us first compute the marginal pdf of \eqref{eq:gspike}.
By integrating out all $x_i$, $i \neq j$, we can easily get:
\bdm
f_{X_j}(x_j)=\frac{1}{\ndim} g(x_j) + \frac{\ndim-1}{\ndim}\delta(x_j).
\edm
Therefore, we have $E[X_j] = 0$.
Now, if $X_i$ and $X_j$ cannot be simultaneously nonzero, therefore,
\bdm
E[X_i X_j] = \delta_{ij} E[X_j^2] = \frac{1}{\ndim}\delta_{ij},
\edm
since the variance of $X_j$ is $1$.
Therefore, the covariance matrix of this process is, as announced,
$I_\ndim/\ndim$.  Therefore, any orthonormal basis is the KLB.
\end{proof}

In other words, the KLB for this process is less restrictive than that
for the simple spike process (Proposition~\ref{prop:spike-KLB}),
and the KLB is again completely useless for this process.

\subsection{Marginal distributions and moments under $\SL^\pm(\ndim,\R)$}
\begin{figure}
\centering{\epsfig{file=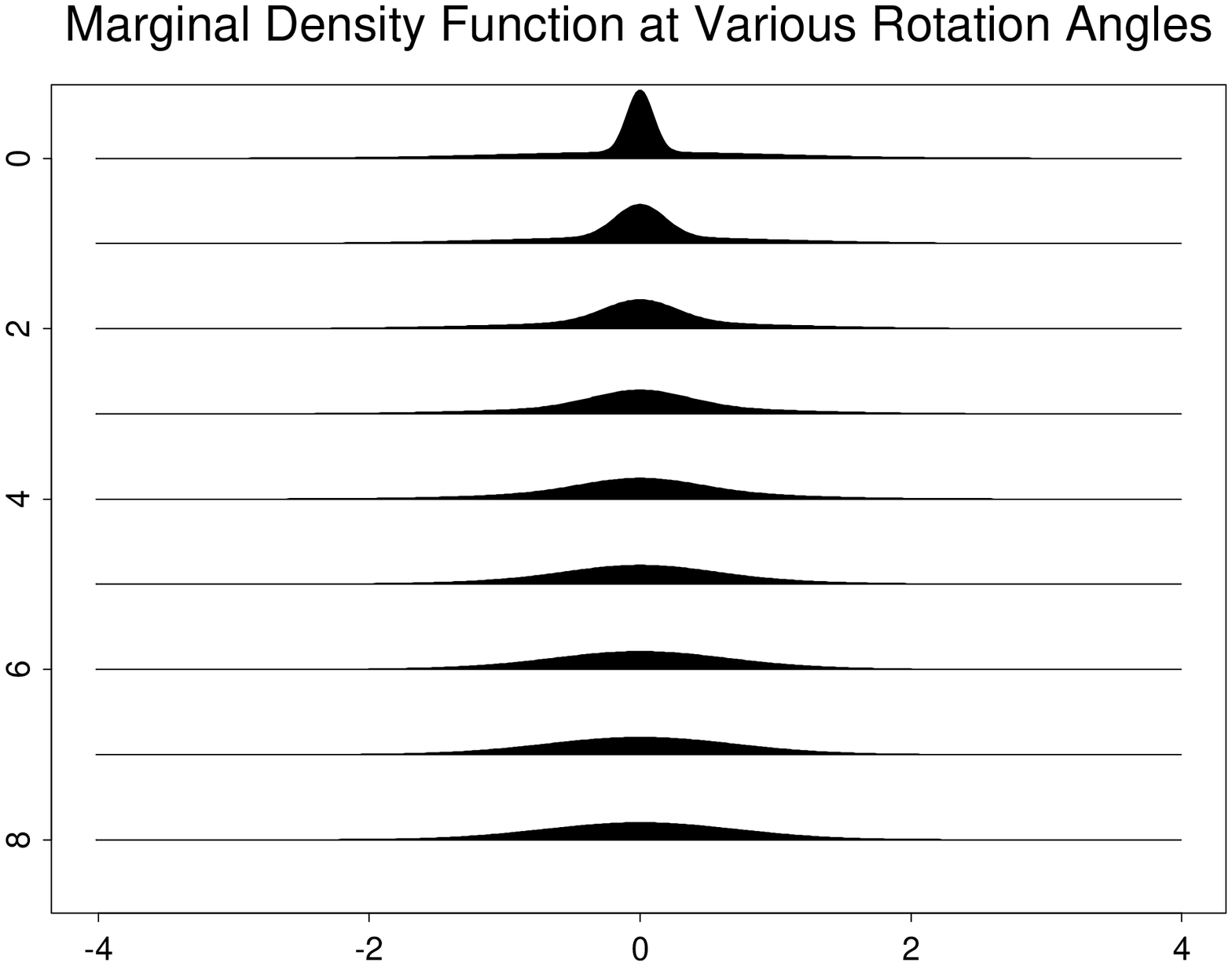,width=0.80\textwidth}}
\caption{The marginal pdf's of the generalized spike process ($\ndim=2$).
All the pdf's shown here are projections of the 2D pdf in Figure~\ref{fig:genspike2} onto the rotated 1D axis. The axis angle in the top row is $0.088$ rad., which is close to the the first axis of the standard basis.  
The axis angle in the bottom row is $\pi/4$ rad., i.e., $45$ degree rotation,
which gives rise to the exact normal distribution.
The other axis angles are equispaced angles between these two.}
\label{fig:mpdf2d}
\end{figure}

\index{$\SL^\pm(\ndim,\R)$}
Before analyzing the BSB and LSDB, we need some background work.
First, let us compute the pdf of the process relative to a transformation
$\bY=B^{-1}\bX$, $B \in \SL^\pm(\ndim,\R)$.
In general, if $\bY=B^{-1}\bX$, then
\bdm
f_\bY(\by) = \frac{1}{|\det B^{-1}|}f_\bX(B\by).
\edm
Therefore, from \eqref{eq:gspike}, and the fact $| \det B |=1$, we have
\begin{equation}
\label{eq:gspike-pdf}
f_\bY(\by) = \frac{1}{\ndim} \sum_{i=1}^\ndim \( \prod_{j \neq i} \delta(\br_j^T\by) \) g(\br_i^T\by),
\end{equation}
where $\br_j^T$ is the $j$th row vector of $B$.
As for its marginal pdf, we have the following lemma:
\index{marginal distribution}
\begin{lem}
\label{lem:gspike-mpdf}
\begin{equation}
\label{eq:gspike-mpdf}
f_{Y_j}(y)= \frac{1}{\ndim} \sum_{i=1}^\ndim g(y; |\Delta_{ij}|),
\quad j=1,\ldots,\ndim,
\end{equation}
where $\Delta_{ij}$ is the $(i,j)$th 
\index{cofactor} cofactor of matrix $B$, and
$g(y;\sigma)=g(y/\sigma)/\sigma$ represents the pdf of the normal distribution
 $\Norm(0,\sigma^2)$. 
\end{lem}
In other words, one can interpret the $j$th marginal pdf as a 
\index{mixture of Gaussians}
\emph{mixture of Gaussians} with the standard deviations $| \Delta_{ij} |$, 
$i=1,\ldots,\ndim$.
Figure~\ref{fig:mpdf2d} shows several marginal pdf's for $\ndim=2$.
As one can see from this figure, it can vary from a very spiky distribution
to a usual normal distribution depending on the rotation angle of 
the coordinate.
\begin{proof}
Let us rewrite \eqref{eq:gspike-pdf} as
\begin{equation}
\label{eq:gspike-pdf2}
f_\bY(\by) = \frac{1}{\ndim} \sum_{i=1}^\ndim 
\delta(\br_1^T \by) \cdots
\delta(\br_{i-1}^T \by) \delta(\br_{i+1}^T \by) \cdots \delta(\br_\ndim^T \by)  
g(\br_i^T\by).
\end{equation}
The $j$th marginal pdf can be written as
\bdm
f_{Y_j}(y_j) = \int f_\bY(y_1,\cdots,y_\ndim) \dd{y_1} \cdots \dd{y_{j-1}} \dd{y_{j+1}} \cdots \dd{y_\ndim}.
\edm
Consider the $i$th term in the summation of \eqref{eq:gspike-pdf2} and 
integrate it out with respect to $y_1,\ldots,y_{j-1},y_{j+1},\ldots,y_\ndim$: 
\begin{equation}
\label{eq:gspike-mpdf-i}
\int \delta(\br_1^T \by) \cdots \delta(\br_{i-1}^T \by) \delta(\br_{i+1}^T \by) \cdots \delta(\br_\ndim^T \by)  g(\br_i^T\by) 
\dd{y_1} \cdots \dd{y_{j-1}} \dd{y_{j+1}} \cdots \dd{y_\ndim}.
\end{equation}
We use a change of variable formula to integrate this.
Let $\br_k^T \by = x_k$, $k=1, \ldots, \ndim$, and let $\bb_\ell$ be
the $\ell$th column vector of $B$.
The relationship $B\by = \bx$ can be rewritten as follows:
\bdm
B^{(i,j)} \by^{(j)} + y_j \bb_j^{(i)} = \bx^{(i)},
\edm
where $B^{(i,j)}$ is the $(\ndim-1) \times (\ndim-1)$ matrix by removing
$i$th row and $j$th column, and the vectors with superscripts indicate 
the length $\ndim-1$ column vectors by removing the elements whose indices
are specified in the parentheses.
This means that
\bdm
\by^{(j)} = \( B^{(i,j)} \)^{-1} \( \bx^{(i)} - y_j \bb_j^{(i)} \).
\edm
Thus,
\begin{eqnarray*}
\dd{\by^{(j)}} & = & \dd{y_1} \cdots \dd{y_{j-1}} \dd{y_{j+1}} \cdots \dd{y_\ndim}\\
&=& \frac{1}{| \det B^{(i,j)} |} \dd{\bx^{(i)}} \\
&=& \frac{1}{| \Delta_{ij} |} \dd{x_1} \cdots \dd{x_{i-1}} \dd{x_{i+1}} \cdots \dd{x_\ndim}.
\end{eqnarray*}
Let us now express $\br_i^T \by=x_i$ in terms of $y_j$ and $\bx$.
\begin{eqnarray}
\label{eq:rTy}
\br_i^T \by &=& \( \br_i^{(j)} \)^T \by^{(j)} + b_{ij} y_j \\ \nonumber
            &=& \( \br_i^{(j)} \)^T  \( B^{(i,j)} \)^{-1} \( \bx^{(i)} - y_j \bb_j^{(i)} \) + b_{ij} y_j \\ \nonumber
            &=& \( \br_i^{(j)} \)^T  \( B^{(i,j)} \)^{-1} \bx^{(i)}
+y_j \( b_{ij} - \( \br_i^{(j)} \)^T  \( B^{(i,j)} \)^{-1} \bb_j^{(i)} \) \\ \nonumber
&\stackrel{(*)}{=}& \( \br_i^{(j)} \)^T  \( B^{(i,j)} \)^{-1} \bx^{(i)}
+ \frac{y_j}{\Delta_{ij}} \det B \\ \nonumber
&=& \( \br_i^{(j)} \)^T  \( B^{(i,j)} \)^{-1} \bx^{(i)} \pm \frac{y_j}{\Delta_{ij}},
\end{eqnarray}
where $(*)$ follows from the following lemma whose proof is shown in
Appendix~\ref{app:detB-delta}:
\begin{lem}
\label{lem:detB-delta}
For any $B=(b_{ij}) \in \GL(\ndim,\R)$,
\bdm
b_{ij} - \( \br_i^{(j)} \)^T  \( B^{(i,j)} \)^{-1} \bb_j^{(i)}
= \frac{1}{\Delta_{ij}} \det B, \quad 1 \leq i, j \leq \ndim.
\edm
\end{lem}
Now, let us go back to the integration \eqref{eq:gspike-mpdf-i}.
Thanks to the property of the delta function with Equation \eqref{eq:rTy},
we have
\begin{eqnarray*}
\int \!\!\!\!\!\! &\cdots& \!\!\!\!\!\! \int \delta(x_1) \cdots \delta(x_{i-1}) \delta(x_{i+1}) \cdots \delta(x_\ndim)  g(\br_i^T\by) \frac{1}{| \Delta_{ij} |} \dd{x_1} \cdots \dd{x_{j-1}} \dd{x_{j+1}} \cdots \dd{x_\ndim}\\
&=& \frac{1}{| \Delta_{ij} |} g( \pm y_j/\Delta_{ij}) \\
&=& g(y_j; |\Delta_{ij}|),
\end{eqnarray*}
where we used the fact that $g(\cdot)$ is an even function.
Therefore, we can write the $j$th marginal distribution as announced in
\eqref{eq:gspike-mpdf}.
\end{proof}

Let us now compute the moments of $Y_i$, which will be used later.
We use the fact that this is \index{mixture of Gaussians} a mixture of $\ndim$ Gaussians each of which
has mean $0$ and variance $| \Delta_{ij} |^2$. 
Therefore, it is obvious to have $E[Y_i]=0$ for all $i=1,\ldots,\ndim$.
Now we have the following lemma for the moments.
\index{moments}
\begin{lem}
\label{lem:p-mom}
\begin{equation}
\label{eq:p-mom}
E[|Y_j|^p] = \frac{\Gamma(p)}{\ndim 2^{p/2-1}\Gamma(p/2)}
\sum_{i=1}^\ndim | \Delta_{ij} |^p, \quad \text{for all $p > 0$.}
\end{equation}
\end{lem}
\begin{proof}
We have:
\begin{eqnarray*}
E[|Y_j|^p] & = & \frac{1}{\ndim} \sum_{i=1}^\ndim \int_{-\infty}^\infty |y|^p g(y;|\Delta_{ij}|) \dd{y} \\
& = & \frac{1}{\ndim} \sum_{i=1}^\ndim \sqrt{\frac{2}{\pi}} |\Delta_{ij}|^p
\Gamma(1+p) D_{-1-p}(0)
\end{eqnarray*}
by Gradshteyn and Ryzhik \cite[Formula 3.462.1]{Gradshteyn-Ryzhik},
where
$D_{-1-p}(\cdot)$ is Whittaker's function as defined by Abramowitz and
Stegun \cite[pp.687]{Abramowitz-Stegun}: 
\bdm
D_{-a-1/2}(0) = U(a,0) = \frac{\sqrt{\pi}}{2^{a/2+1/4} \, \Gamma(a/2+3/4)}.
\edm
Thus, putting $a=p+1/2$ to the above equation yields:
\bdm
D_{-1-p}(0) = \frac{\sqrt{\pi}}{2^{1/2+p/2} \, \Gamma(1+p/2)}.
\edm

Therefore, we have
\begin{eqnarray*}
E[|Y_j|^p] & = & \frac{1}{\ndim} \sum_{i=1}^\ndim | \Delta_{ij} |^p
\frac{\Gamma(1+p)}{2^{p/2} \, \Gamma(1+p/2)} \\
& = & \frac{1}{\ndim} \sum_{i=1}^\ndim | \Delta_{ij} |^p
\frac{\Gamma(p)}{2^{p/2-1} \, \Gamma(p/2)} \\
& = & \frac{\Gamma(p)}{\ndim 2^{p/2-1} \, \Gamma(p/2)}
\sum_{i=1}^\ndim | \Delta_{ij} |^p,
\end{eqnarray*}
as we desired.
\end{proof}

\subsection{The Best Sparsifying Basis}
\index{best sparsifying basis (BSB)}
As for the BSB, after all, there is no difference between the generalized 
spike process and the simple spike process.
\begin{thm}
\label{thm:gspike-sparse}
The BSB with any $p \in [0,1]$ for the generalized spike process is
the \index{standard basis} standard basis 
if \index{$\OG(\ndim)$} $\Dict=\OG(\ndim)$ or 
\index{$\SL^\pm(\ndim,\R)$} $\SL^\pm(\ndim,\R)$.
\end{thm}
\begin{proof}
Let us first consider the case $p \in (0,1]$.
Then, using Lemma~\ref{lem:p-mom}, the cost function \eqref{eq:cp} can be
rewritten as follows:
\bdm
\Cost_p(B \cond \bx) = \sum_{j=1}^\ndim E[ | Y_j |^p ]
= \frac{\Gamma(p)}{\ndim 2^{p/2-1} \, \Gamma(p/2)}
\sum_{i=1}^\ndim \sum_{j=1}^\ndim | \Delta_{ij} |^p.
\edm
Let us now define a matrix $\tilde{B} \define ( \Delta_{ij} )$.
Then $\tilde{B} \in \SL^\pm(\ndim,\R)$ since 
\bdm
B^{-1}=\frac{1}{\det B} ( \Delta_{ji} ) = \pm ( \Delta_{ji} ),
\edm
and $B^{-1} \in \SL^\pm(\ndim,\R)$.
Therefore, this reduces to
\bdm
\Cost_p(B \cond \bx) = \frac{\Gamma(p)}{\ndim 2^{p/2-1} \, \Gamma(p/2)}
\sum_{i=1}^\ndim \sum_{j=1}^\ndim | \tilde{b}_{ij} |^p 
= \Cost_p(\tilde{B} \cond \bx).
\edm
This means that our problem now becomes the same as Theorem~1 
in \cite{BENICHOU-SAITO-APBOOK} (or Theorem~\ref{thm:spike-sparse} in
this paper) by replacing $B$ by $\tilde{B}$.
Thus, it asserts that the $\tilde{B}$ must be the identity matrix $I_\ndim$
or its permuted or sign flipped versions.
Suppose $\Delta_{ij}=\delta_{ij}$.
Then, $B^{-1}=\pm(\Delta_{ji})=\pm I_\ndim$, which implies that
$B=\pm I_\ndim$.  If $( \Delta_{ji} )$ is any permutation matrix, then
$B^{-1}$ is just that permutation matrix or its sign flipped version.
Therefore, $B$ is also a permutation matrix or its sign flipped version.

Finally, let us consider the case $p=0$.
Then, any linear invertible transformation except the identity matrix
or its permuted or sign-flipped versions clearly increases the number of
nonzero elements after the transformation.
Therefore, the BSB with $p=0$ is also a permutation matrix or its sign flipped
version.

This completes the proof of Theorem~\ref{thm:gspike-sparse}.
\end{proof}

\subsection{The LSDB/KMB among $\OG(\ndim)$}
\index{kurtosis-maximizing basis (KMB)}
\index{least statistically-dependent basis (LSDB)}
\index{$\OG(\ndim)$}
As for the LSDB/KMB, we can see some difference from the simple spike process.

Let us now consider a more specific case of $\Dict=\OG(\ndim)$.
So far, we have been unable to prove the following conjecture.
\begin{conj}
\label{conj:gspike-OG}
The LSDB among $\OG(\ndim)$ is the standard basis.
\end{conj}
The difficulty is the evaluation of the sum of the marginal entropies
\eqref{eq:c-lsdb} for the pdf's of the form \eqref{eq:gspike-mpdf}.
However, a major simplification occurs if we consider the KMB instead of 
the LSDB, and we can prove the following:
\begin{thm}
\label{thm:gspike-kmb-OG}
The KMB among $\OG(\ndim)$ is the \index{standard basis} standard basis.
\end{thm}
\begin{proof}
Because $E[Y_j]=0$ and $E[Y_j^2]=\frac{1}{\ndim}\sum_{i=1}^\ndim \Delta_{ij}^2$
for all $j$, the fourth order central moment of $Y_j$ can be written as
$\mu_4(Y_j) = \frac{3}{\ndim} \sum_{i=1}^\ndim \Delta_{ij}^4$,
and consequently the cost function in \eqref{eq:c-kmb} becomes
\begin{equation}
\label{eq:gspike-c-kmb-OG}
\Cost_\kappa(B \cond \bX) = \frac{3}{\ndim} \sum_{j=1}^\ndim \(  \sum_{i=1}^\ndim \Delta_{ij}^4 - \frac{1}{\ndim} \(\sum_{i=1}^\ndim \Delta_{ij}^2 \)^2 \).
\end{equation}
Note that this is true for any $B \in \SL^\pm(\ndim,\R)$.
If we restrict our basis search within $\OG(\ndim)$, another major 
simplification occurs because we have the following special relationship 
between $\Delta_{ij}$ and the matrix element $b_{ji}$ of $B \in \OG(\ndim)$:
\bdm
B^{-1}=\frac{1}{\det B} \( \Delta_{ji} \) = B^T.
\edm
In other words,
\bdm
\Delta_{ij} = (\det B) b_{ij} = \pm b_{ij}.
\edm
Therefore, we have
\bdm
\sum_{i=1}^\ndim \Delta_{ij}^2=\sum_{i=1}^\ndim b_{ij}^2=1.
\edm
Inserting this into \eqref{eq:gspike-c-kmb-OG}, we get the following
simplified cost for $\Dict=\OG(\ndim)$:
\bdm
\Cost_\kappa(B \cond \bX) =
-\frac{3}{\ndim} \(1 - \sum_{i=1}^\ndim \sum_{j=1}^\ndim \Delta_{ij}^4 \). 
\edm
This means that the KMB can be rewritten as follows:
\begin{equation}
\label{eq:gspike-kmb-OG}
B_\kappa = \arg \max_{B \in \OG(\ndim)} \sum_{i,j} b_{ij}^4.
\end{equation}
Let us note that the existence of the maximum is guaranteed
because the set $\OG(\ndim)$ is \emph{compact} and
the cost function $\sum_{i,j} b_{ij}^4$ is continuous, 

Now, let us consider a matrix $P=(p_{ij})=(b^2_{ij})$.
Then, from the orthonormality of columns and rows of $B$, this matrix $P$
belongs to a set of \index{doubly stochastic matrices}
\emph{doubly stochastic matrices} $\DS(\ndim)$.
Since doubly stochastic matrices obtained by squaring the elements 
of $\OG(\ndim)$ consist of a proper subset of $\DS(\ndim)$, we have
\bdm
\max_{B \in \OG(\ndim)} \sum_{i,j} b_{ij}^4 \leq \max_{P \in \DS(\ndim)} \sum_{i,j} p_{ij}^2.
\edm
Now, we prove that such $P$ must be an identity matrix or its permuted 
version.
\begin{eqnarray*}
\max_{P \in \DS(\ndim)} \sum_{j=1}^\ndim \sum_{i=1}^\ndim  p_{ij}^2 &\leq& \sum_{j=1}^\ndim 
\( \max_{\sum_{i=1}^\ndim p_{ij}=1} \sum_{i=1}^\ndim p_{ij}^2 \) \\
&=& \sum_{j=1}^\ndim 1\\
&=& \ndim,
\end{eqnarray*}
where the first equality follows from the fact that maxima of the
radius of the sphere $\sum_i p_{ij}^2$ subject to
$\sum_i p_{ij}=1$, $p_{ij} \geq 0$ occur only at the vertices of
that simplex, i.e., $\bp_j=\be_{\sigma(j)}$, $j=1,\ldots,\ndim$ where
$\sigma(\cdot)$ is a permutation of $\ndim$ items.
That is, the column vectors of $P$ must be the standard basis vectors.
This implies that the matrix $B$ corresponding to $P=I_\ndim$ or its
permuted version must be either $I_\ndim$ or its permuted and/or sign-flipped 
version.
\end{proof}

\subsection{The LSDB/KMB among $\SL^\pm(\ndim,\R)$}
\index{kurtosis-maximizing basis (KMB)}
\index{least statistically-dependent basis (LSDB)}
\index{$\SL^\pm(\ndim,\R)$}
If we extend our search to this more general case, 
we have the following theorem.
\begin{thm}
\label{thm:gspike-kmb-SL}
The KMB among $\SL^\pm(\ndim,\R)$ does not exist.
\end{thm}
\begin{proof}
The set $\SL^\pm(\ndim,\R)$ is not compact.  
Therefore, there is no guarantee that the cost function 
$\Cost_\kappa(B \cond \bX)$ has a minimum value on this set.
One can in fact consider a simple counter-example,
$B=\diag(a,a^{-1},1,\cdots,1)$, where $a$ is any nonzero real scalar.
Then, one can show that $\Cost_\kappa(B \cond \bX)=-(a^4+a^{-4}+\ndim-2)$,
which tends to $-\infty$ as $a \uparrow \infty$.
\end{proof}

As for the LSDB, we do not know whether the LSDB exists among
$\SL^\pm(\ndim,\R)$ at this point, although
we believe that the LSDB is the standard basis (or its permuted/sign-flipped
versions).
The negative result in the KMB does not imply the negative result in
the LSDB.

\section{Discussion}
\label{sec:disc}
Unlike the simple spike process, \index{best sparsifying basis (BSB)}
the BSB and \index{kurtosis-maximizing basis (KMB)} the KMB (an alternative
to the LSDB) selects \index{standard basis} the standard basis if we restrict
our basis search within $\OG(\ndim)$.
If we extend our basis search to \index{$\SL^\pm(\ndim,\R)$} $\SL^\pm(\ndim,\R)$, then the BSB exists and is again the standard basis whereas
the KMB does not exist.

Although the generalized spike process is a simple stochastic process,
we have the following important interpretation.
Consider a stochastic process generating a basis vector randomly selected from
some fixed orthonormal basis and multiplied by a scalar varying as the 
standard normal distribution at a time.
Then, both that basis itself is the BSB and the KMB among $\OG(\ndim)$.
Theorems~\ref{thm:gspike-sparse} and \ref{thm:gspike-kmb-OG} claim that 
once we transform the data to the generalized spikes, one cannot do any 
better than that both in sparsity and independence within $\OG(\ndim)$.
Of course, if one extends the search to nonlinear transformations,
then it becomes a different story.  We refer the reader to our
recent articles \cite{LIN-SAITO-LEVINE-INGA}, \cite{LIN-SAITO-LEVINE-INGA2}, 
for the details of a nonlinear algorithm.

The results of this paper further support our conclusion of the
previous paper: dealing with the BSB is much simpler than the LSDB.
To deal with statistical dependency, we need to consider the probability law
of the underlying process (e.g., entropy or the marginal pdf's) explicitly. 
That is why we need to consider the KMB instead of the LSDB to prove the
theorems.
Also in practice, given a finite set of training data, it is a nontrivial task
to reliably estimate the marginal pdf's.
Moreover, the LSDB unfortunately cannot tell how close it is to the true 
statistical independence; it can only tell that it is the best one (i.e., 
the closest one to the statistical independence) among the given set of 
possible bases.
In order to quantify the absolute statistical dependence,
we need to estimate the true high-dimensional entropy of the original process,
$H(\bX)$, which is an extremely difficult task in general.
We would like to note, however, a recent attempt to estimate the 
high-dimensional entropy of the process by Hero and Michel \cite{HERO-MICHEL},
which uses the minimum spanning trees of the input data and does not
require to estimate the pdf of the process.
We feel that this type of techniques will help assessing the absolute
statistical dependence of the process under the LSDB coordinates.
Another interesting observation is that the KMB is rather sensitive to
the \index{orthonormality} orthonormality of the basis dictionary whereas the BSB is insensitive
to that.  Our previous results on the simple spike process (e.g.,
Theorems~\ref{thm:spike-OG}, \ref{thm:spike-GL}) also suggest the sensitivity
of the LSDB to the orthonormality of the basis dictionary.
This may restrict and discourage us to develop a new basis or a new basis 
dictionary that optimize the statistical independence.

On the other hand, the sparsity criterion neither requires estimating
the marginal pdf's nor reveals the sensitivity to the orthonormality.
Simply computing the expected $\ell^p$ norms suffices.
Moreover, one can even adapt the BSB for each realization rather than for
the whole realizations, which is impossible for the LSDB, as we discussed in
\cite{BENICHOU-SAITO-APBOOK}, \cite{SAITO-LARSON-BENICHOU}, \cite{SAITO-BENICHOU-LARSON-LEBORNE-LUCERO}.

These observations, therefore, suggest that the pursuit of sparse 
representations should be encouraged rather than that of statistically 
independent representations, if we believe that mammalian vision systems 
were evolved and developed by the principle of data compression.
This is also the view point indicated by Donoho \cite{DONOHO-SCA}.

Finally, there are a few interesting generalizations of the spike processes,
which need to be addressed in the near future.
We need to consider a stochastic process that randomly throws in multiple 
spikes to a single realization.
If one throws in more and more spikes to one realization, 
the \index{standard basis} standard basis
is getting worse in terms of sparsity.
Also, we can consider various rules to throw in multiple spikes.
For example, for each realization, we can select the locations of the spikes
statistically independently.  This is the simplest multiple spike process.
Alternatively, we can consider a certain dependence in choosing the locations
of the spikes.  The ramp process of Yves Meyer analyzed by the wavelet basis
is such an example; each realization of the ramp process generates a small
number of spikes in the wavelet coefficients in the locations determined by
the location of the discontinuity of the process.
See \cite{BUCK-DONO}, \cite{DVDD}, \cite{Meyer-NewBook}, \cite{SAITO-LARSON-BENICHOU} for more about the ramp process.

Unless very special circumstances, it would be extremely difficult to
find the BSB of a complicated stochastic process (e.g., natural scene
images) that truly converts its realizations to the spike process.
More likely, a theoretically and computationally feasible basis that 
sparsifies the realizations of a complicated process well (e.g., curvelets
for the natural scene images \cite{DONOHO-FLESIA}) may generate
expansion coefficients that may be viewed as 
\index{multiple spike process}
an amplitude-varying multiple spike process.
In order to tackle this scenario, we certainly need to: 1) identify 
interesting, useful, and simple enough specific stochastic processes;
2) develop the BSB adapted to such specific processes; and 3) deepen our 
understanding of the amplitude-varying multiple spike process.

\section*{Acknowledgment}
I would like to thank the fruitful discussions with Dr.\ Motohico Mulase
and Dr.\ Roger Wets, of UC Davis.
This research was partially supported by NSF DMS-99-73032, DMS-99-78321,
and ONR YIP N00014-00-1-046.

\appendix
\section{Proof of Lemma~\ref{lem:detB-delta}}
\label{app:detB-delta}
\begin{proof}
Let us consider the following system of linear equations:
\bdm
B^{(i,j)} \bz^{(j)} = \bb_j^{(i)},
\edm
where $\bz^{(j)}=(z_1,\cdots,z_{j-1},z_{j+1},\cdots,z_\ndim)^T \in \R^{\ndim-1}$,
$j=1,\ldots,\ndim$.
Using Cramer's rule (e.g., \cite[pp.21]{HORN-JOHNSON}),
we have, for $k=1,\ldots,j-1,j+1,\ldots,\ndim$,
\begin{eqnarray*}
z^{(j)}_k &=& \frac{1}{\det B^{(i,j)}} \det \left[
\bb_1^{(i)} \matpart \cdots \matpart \bb_{k-1}^{(i)} \matpart \bb_j^{(i)} \matpart \bb_{k+1}^{(i)} \matpart \cdots \matpart \bb_\ndim^{(i)} \right] \\
&\stackrel{(a)}{=}& (-1)^{|k-j|-1} \frac{B^{(i,k)}}{B^{(i,j)}}\\
&\stackrel{(b)}{=}& (-1)^{|k-j|-1} \frac{\Delta_{ik}/(-1)^{i+k}}{\Delta_{ij}/(-1)^{i+j}}\\
&=& - \frac{\Delta_{ik}}{\Delta_{ij}},
\end{eqnarray*}
where $(a)$ follows from the $(|k-j|-1)$ column 
permutations to move $\bb_j^{(i)}$ located at the $k$th column to the $j$th
column of $B^{(i,j)}$, and $(b)$ follows from the definition of the cofactor.
Hence, 
\begin{eqnarray*}
b_{ij} - \( \br_i^{(j)} \)^T  \( B^{(i,j)} \)^{-1} \bb_j^{(i)}
&=& b_{ij} - \( \br_i^{(j)} \)^T  \bz^{(j)}\\
&=& b_{ij} + \frac{1}{\Delta_{ij}} \sum_{k \neq j} b_{ik} \Delta_{ik}\\
&=& \frac{1}{\Delta_{ij}} \sum_{k=1}^\ndim b_{ik} \Delta_{ik}\\
&=& \frac{1}{\Delta_{ij}} \det B.
\end{eqnarray*}
This completes the proof of Lemma~\ref{lem:detB-delta}.
\end{proof}

\bibliographystyle{plain}

\end{document}